\documentclass[12pt,a4paper]{amsart}
      \usepackage[english]{babel}
      \usepackage{amsmath}
      \usepackage{amstext}
      \usepackage{amssymb}
      \usepackage{amsthm}
      \usepackage{amsfonts}
      \usepackage[T1]{fontenc}
      \usepackage{latexsym}
      \usepackage{psfrag}
      \newtheorem{thm}{Theorem}[section]
      \newtheorem{lem}[thm]{Lemma}
      \newtheorem{prop}[thm]{Proposition}

      \newtheorem{coro}[thm]{Corollary}
      
      \newtheorem{example}[thm]{Example}

      \newcommand{\N}{\mathbb N}
      \newcommand{\R}{\mathbb R}
      
      \newcommand{\Z}{\mathbb Z}
      \newcommand{\ds}{\displaystyle}
      \newcommand{\U}{\mathbb U}
      
      \newcommand{\Q}{\mathbb Q}

      \newenvironment{appli}{\left( \begin{array}{ccc}}{\end{array} \right)}

      \newcommand{\ba}{\begin{appli}}
      \newcommand{\ea}{\end{appli}}

\begin{document}
\title{Geometry in Urysohn's universal metric space}
\author{Julien Melleray}
\thanks{MSC: 51F99}
\date{}
\maketitle

\begin{abstract} \noindent In recent years, much interest was
devoted to the Urysohn space and its isometry group; this paper is
a contribution to this field of research. We mostly concern
ourselves with the properties of isometries of $\U$, showing for
instance that any Polish metric space is isometric to the set of
fixed points of some isometry $\varphi$. We conclude the paper by
studying a question of Urysohn, proving that compact homogeneity
is the strongest homogeneity property possible in $\U$.
\end{abstract}

\begin{section}{Introduction} \label{intro}
\noindent In a paper published posthumously (see \cite{Urysohn}),
P.S Urysohn constructed a complete separable metric space $\U$
that is \textit{universal}, i.e contains an isometric copy of
every complete separable metric space. This seems to have been
forgotten for a while, perhaps because around
the same time Banach and Mazur proved that ${\mathcal C}([0,1])$ is also universal. \\
Yet, the interest of the Urysohn space $\U$ does not lie in its
universality alone: as Urysohn himself had remarked, $\U$ is also
$\omega$-\textit{homogeneous}, i.e for any two finite subsets $A$,
$B$ of $\U$ which are isometric (as abstract metric spaces), there
exists an isometry $\varphi$ of $\U$ such that $\varphi(A)=B$.
Moreover, Urysohn proved that $\U$ is, up to isometry, the only
universal
$\omega$-homogeneous Polish metric space.\\
In the case of Polish metric spaces, it turns out that
universality and $\omega$-homogeneity can be merged in one
property, called \textit{finite injectivity}: a metric space
$(X,d)$ is finitely injective iff for any pair of finite metric
spaces $K\subseteq L $ and any isometric embedding $\varphi \colon
K \to X$, there exists an isometric embedding $\tilde \varphi
\colon L \to X$ such that $\tilde {\varphi}_{|_K}=\varphi$.\\
Then one can prove that a Polish metric space is universal and
$\omega$-homogeneous if, and only if, it is finitely injective;
this is also due to Urysohn, who was the first to use finite
injectivity (using another definition of it). \footnote{About
finite injectivity, Urysohn stated in \cite{Urysohn} "Voici la
propriété fondamentale de l'espace $\U$ dont, malgré son caractère
auxiliaire, les autres
propriétés de cet espace sont des conséquences plus ou moins immédiates".} \\
This point of view highlights the parallel between $\U$ and other
universal objects, such as the universal graph for instance; the
interested reader can find a more detailed exposition of this and references in \cite{camver}.\\
The interest in $\U$ was revived in 1986 when Kat$\check
{\mbox{e}}$tov, while working on analogues of the Urysohn space
for metric spaces of a given density character,  gave in
\cite{Katetov} a new construction of $\U$, which enables one to
naturally "build" an isometric copy of $\U$ "around" any separable
metric space $X$. In \cite{Uspenskij2} Uspenskij remarked that
this construction (which we will detail a bit more in section
\ref{nota}) enables one to keep track of the isometries of $X$,
and used that to obtain a canonical continuous embedding of the
group of isometries of $X$ into $Iso(\U)$, the group of isometries
of $\U$ (both groups being endowed with the product topology,
which turns $Iso(\U)$ into a Polish group). Since any Polish group
$G$ continuously embeds in the isometry group of some Polish space
$X$, this shows that any Polish group
is isomorphic to a (necessarily closed) subgroup of $Iso(\U)$.\\
This result spurred interest for the study of $\U$; in
\cite{Vershik}, Vershik showed that generically (for a natural
Polish topology on the sets of distances on $\N$) the completion
of a countable metric space is finitely injective, and thus
isometric to $\U$; in \cite{Uspenskij3} Uspenskij completely
characterized the topology of $\U$ by
showing, using Torunczyk's criterion, that $\U$ is homeomorphic to $l^2(\N)$.\\
During the same period, Gao and Kechris used $\U$ to study the
complexity of the equivalence relation of isometry between certain
classes of Polish metric spaces (viewed as elements of $\mathcal
{F}(\U)$). For instance, they proved that the relation of isometry
between Polish metric spaces is Borel bi-reducible to the
translation action of $Iso(\U)$ on ${\mathcal F}(\U)$, given by
$\varphi.F=\varphi(F)$, and that this relation is universal among
relations induced by a continuous action of a Polish group (see
\cite{gaokec} for a detailed exposition of their results
and references about the theory of Borel equivalence relations).\\
Despite all the recent interest in $\U$, not much work has yet
been done on its geometric properties, with the exception of
\cite{camver}, where the authors build interesting examples of
subgroups of $Iso(\U)$. \\
As Urysohn himself had understood, finite injectivity has
remarkable consequences on the geometry of $\U$, some of which we
study in section \ref{geom}; we begin with the easy fact that any
isometry map which coincides with $id_{\U}$ on a set of non-empty
interior must actually be $id_{\U}$. We then go on to study a bit
the isometric copies of $\U$ contained in $\U$, e.g we show that
$\U$ is isometric to $\U \setminus B$, where $B$ is any open ball in $\U$. \\
We also use similar ideas to study the sets of fixed points of
isometries, proving in particular that any Polish metric space is
isometric to the set of fixed points of some isometry of
$\U$.\\
The remainder of the article is devoted to the study of a question
of Urysohn, who asked in \cite{Urysohn} whether $\U$ had stronger
homogeneity properties than $\omega$-homogeneity \footnote{"On
demandera, peut-être, si l'espace $\U$ ne jouit pas d'une
propriété d'homogénéité plus précise que celle que nous avons
indiquée au n. 14". }; we build on known results to solve that
problem. Most importantly, we use the tools introduced by
Kat\v{e}tov in \cite{Katetov}.
Let us state precisely the problems we concern ourselves with:\\

\noindent {\bf Question $\mathbf 1$.} Characterize the Polish
metric spaces $(X,d)$ such that whenever $X_1,X_2 \subseteq \U$
are isometric to $X$, there is an isometry $\varphi$ of $\U$ such
that $\varphi(X_1)=X_2$.\\

\noindent As it turns out, we will not directly study that
question, but another related one, which can be thought of as
looking if one can extend finite injectivity:\\

\noindent {\bf Question 2} Characterize the Polish metric spaces
$(X,d)$ such that, whenever $X' \subseteq \U$ is isometric to $X$
and $f \in E(X')$, there is $z \in \U$ such that $\forall x \in
X'$, $d(x,z)=f(x)$.\\
($E(X)$ denotes the set of Kat\v{e}tov maps on $X$).\\

\noindent It is rather simple, as we will see in section
\ref{trans}, to show that Property 1 implies Property 2, and it is
a well-known fact (see \cite{huhu} or \cite{Gromov}) that the
answer to both questions is positive whenever $X$ is compact:
\begin{thm} (Huhunai\v{s}vili) \label{huhu}
If $K \subseteq \U$ is compact and $f \in E(K)$, then there is $z
\in \U$ such that $\ d(z,x)=f(x)$ for all $x\in K$.
\end{thm}
\begin{coro}
If $K,L \subseteq \U$ are compact and $\varphi \colon K \to L$ is
an isometry, then there is an isometry $\tilde \varphi \colon \U
\to \U$ such that $\tilde \varphi_{|_K}=\varphi$.
\end{coro}
\noindent The corollary is deduced from the theorem by the
standard back-and-forth method (So, in that case, a positive answer to question 2
enables one to answer positively question 1; we will see that it is actually always the case).\\
Remarking that if $X$ is such that $E(X)$ is not separable then
$X$ can have neither property (1) nor property (2), we provide a
characterization of the spaces $X$ such that $E(X)$ is separable,
which we tentatively call \textit{compactly tentacular spaces},
for reasons that should be explained in a later version of the
paper. Afterwards, we show that, if $X$ is not compact
and is compactly tentacular then $X$ does not have property 2 either.\\
Therefore, our results enable us to deduce that a space has
property 1 (or 2) if, and only if, it is compact, thus answering
Urysohn's question: compact homogeneity is the strongest
homogeneity property possible in $\U$.\\

\noindent \textit{Acknowledgements.} I would not have written this
article if not for many conversations with Thierry Monteil; he
introduced me to the Urysohn space, and the results in section
\ref{geom} are answers to questions we asked during these
conversations. I have learnt much from these talks, and for that I
am extremely grateful; the results in section \ref{geom} are partly his.\\
I also would like to seize the opportunity to thank Vladimir
Pestov, whose kindness, disponibility and insights while at the
CIRM in september 2004 were much appreciated; and Alekos Kechris,
for numerous remarks and suggestions. Last (but not least of
course), I would like to thank all those at the Equipe d'Analyse
Fonctionnelle of the Université Paris 6, especially the members of
the descriptive set theory workgroup, for providing me with such a
nice environment to work in. \\

\end{section}


\begin{section}{Notations and definitions} \label{nota}
\noindent If $(X,d)$ is a complete separable metric space, we say
that it is a \textit{Polish metric space}, and often
write it simply $X$. \\
If $X$ is a topological space and there is a distance $d$ on $X$
which induces the topology of $X$ and is such
that $(X,d)$ is a Polish metric space, we say that the topology of $X$ is Polish.  \vspace{0.2cm}\\
If $(X,d)$ is a metric space, $x \in X$ and $r >0$, we use the
notation $B(x,r[$ (resp. $B(x,r]$ ) to denote the open
(resp. closed) ball of center $x$ and radius $r$; $S(x,r)$ denotes the sphere of center $x$ and radius $r$.\vspace{0.2cm}\\
To avoid confusions, we say, if $(X,d)$ and $(X',d')$ are two
metric spaces and $f$ is a map from $X$ into $X'$, that $f$ is an
\textit{isometric map} if $d(x,y)=d'(f(x),f(y))$ for all $x,y \in
X$. If additionally $f$ is onto, then we say that $f$ is an
\textit{isometry}. \vspace{0.2cm}\\
A \textit{Polish group} is a topological group whose topology is
Polish; if $X$ is a separable metric space, then we denote its
isometry group by $Iso(X)$, and endow it with the product
topology, which turns it into a second countable topological
group, and into a Polish group if $X$ is Polish (see \cite{beckec}
or \cite{Kechris1} for a thorough introduction to the theory of Polish groups). \vspace{0.2cm}\\
\noindent Let $(X,d)$ be a metric space; we say that $f: X \to \R
$ is a \textit{Kat\v{e}tov map } if
$$\forall x,y \in X \ |f(x)-f(y)|\leq d(x,y) \leq f(x)+f(y)\ \ .$$
These maps correspond to one-point metric extensions of $X$. We
denote by $E(X)$ the set of all Kat\v{e}tov maps on $X$; we endow
it with the sup-metric, which turns it into a complete
metric space.\\
That definition was introduced by Kat$\check {\mbox{e}}$tov in
\cite{Katetov}, and it turns out to be pertinent to the study of
finitely injective spaces, since one can easily see by induction
that a non-empty metric space $X$ is finitely injective if, and
only if,
$$\forall x_1,\ldots,x_n \in X \ \forall f \in E(\{x_1,\ldots,x_n\})\ \exists z \in X\ \forall x \in X \ d(z,x)=f(x)\ .$$
(This is the form under which Urysohn used finite injectivity in
his original article). \\
If $Y \subseteq X$ and $f \in E(Y)$, define ${k}(f)\colon X \to
\R$ ( the Kat\v{e}tov extension of $f$) by
$${k}(f)(x)= \inf\{f(y)+d(x,y) \colon y \in Y\} .$$
Then ${k}(f)$ is the greatest 1-Lipschitz map on $X$ which is
equal to $f$ on $Y$; one checks easily (see for instance
\cite{Katetov}) that ${k}(f) \in E(X) $ and $f\mapsto {k}(f)$ is
an isometric embedding of $E(Y)$ into $E(X)$. \vspace{0.2cm}\\
To simplify future definitions, if $f \in E(X)$ and $S \subseteq X$ are such that \\
$ \ f(x)=\inf\{f(s)+d(x,s) \colon s \in S\}$ for all $x \in X$, we
say that $S$ is a \textit{support of }$f$, or that $S$ \textit{controls} $f$. \\
Notice that if $S$ controls $f\in E(X)$ and $S \subseteq T$, then $T$ controls $f$. \vspace{0.2cm}\\
Similarly, $X$ isometrically embeds in $E(X)$ via the Kuratowski
map $x \mapsto f_x$, where $f_x(y)=d(x,y)$. A crucial fact for our
purposes is that
 $$ \forall f \in E(X)\ \forall x \in X\  d(f,f_x)=f(x).$$
Thus, if one identifies $X$ to a subset of $E(X)$ via the
Kuratowski map, $E(X)$ is a metric space containing $X$ and such
that all one-point metric extensions of $X$ embed isometrically in
$E(X)$.\\

\noindent We now go on to sketching Kat\v{e}tov's construction of
$\U$; we refer the reader to \cite{gaokec}, \cite{Gromov},
\cite{Katetov} or \cite{Uspenskij2} for a more detailed
presentation
and proofs of the results we will use below.\\
Most important for the construction is the following
\begin{thm}\label{completion}(Urysohn)
If $X$ is a finitely injective metric space, then the completion of $X$
is also finitely injective.
\end{thm}
\noindent Since $\U$ is, up to isometry, the unique finitely injective Polish metric space,
this proves that the completion of any separable finitely injective metric space is isometric to
$\U$.\\
The basic idea of Kat\v{e}tov's construction works like this: if
one lets $X_0=X$, $X_{i+1}=E(X_i)$ then, identifying each $X_i$ to
a subset of $X_{i+1}$ via the Kuratowski map,
let $Y$ be the inductive limit of the sequence $X_i$. \\
The definition of $Y$ makes it clear that $Y$ is finitely
injective, since any $\{x_1,\ldots,x_n\}\subseteq Y$ must be
contained in some $X_m$, so that for any $f\in
E(\{x_1,\ldots,x_n\})$ there exists $z \in X_{m+1}$ such that
$d(z,x_i) =f(x_i)$ for all $i$. \\
Thus, if $Y$ were separable, its completion would be isometric to
$\U$, and
one would have obtained an isometric embedding of $X$ into $\U$.\\
The problem is that $E(X)$ is in general not separable (see
section \ref{trans}).\\
At each step, we have added too many functions; define then
$$E(X,\omega)=\{f \in E(X)\colon f \mbox{ is controlled by some finite } S\subseteq X\}\ .$$
Then $E(X,\omega)$ is easily seen to be separable if $X$ is, and the Kuratowski map actually maps $X$
into $E(X,\omega)$, since each $f_x$ is controlled by $\{x\}$.\\
Notice also that, if $\{x_1,\ldots,x_n\} \subseteq X$ and $f \in
E(\{x_1,\ldots,x_n\})$, then its Kat\v{e}tov extension $ {k}(f)$
is in $E(X,\omega)$, and $d({k}(f),f_{x_i})=f(x_i)$ for all
$i$. \vspace{0.2cm}\\
Thus, if one defines this time $X_0=X$, $X_{i+1}=E(X_i,\omega)$,
and assume again that $X_i \subseteq X_{i+1}$ then $Y=\cup X_i$ is
separable and finitely injective, hence its completion $Z$ is
isometric to $\U$, and $X \subseteq Z$. \\
The most interesting property of this construction is that it
enables one to keep track of the isometries of $X$: indeed, any
$\varphi\in Iso(X)$ is the restriction of a unique isometry
$\tilde \varphi$ of $E(X,\omega)$, and the mapping $\varphi
\mapsto \tilde \varphi$ from $Iso(X)$ into $Iso(E(X,\omega))$ is a
continuous group embedding (see \cite{Katetov}).\\
That way, we obtain for all $i \in \N$ continuous embeddings
$\Psi^i\colon Iso(X) \to Iso(X_i)$, such that
$\Psi^{i+1}(\varphi)_{|_{X_i}}=\Psi^i(\varphi)$ for all $i$
and all $\varphi \in Iso(X)$.\\
This in turns defines a continuous embedding from $Iso(X)$ into
$Iso(Y)$, and since extension of isometries defines a continuous
embedding from the isometry group of any metric space into that of
its completion (see \cite{Uspenskij1}), we actually have a
continuous embedding of $Iso(X)$ into the isometry group of $Z$,
that is to say $Iso(\U)$
(and the image of any $\varphi \in Iso(X)$ is actually an extension of $\varphi$ to $\U$ ).\vspace{0.2cm}\\
In the remainder of the text, we follow \cite{pest} and say that a
metric space $X$ is \textit{g-embedded} in $\U$ if $X$ is embedded
in $\U$, and there is a continuous morphism $\Phi \colon Iso(X)
\to Iso(\U)$ such that $\Phi(\varphi)$
extends $\varphi$ for all $\varphi \in Iso(X)$.\\
\end{section}


\begin{section}{Finite injectivity and the geometry of
$\U$}\label{geom}
\begin{subsection}{First results}$ $\\
The following result, tough easy to prove, is worth stating on its
own, since it gives a good idea of the kind of problems we concern
ourselves with in this section:
\begin{thm} \label{regul}
  If $\varphi\colon \U \to \U$ is an isometric map, and $\varphi_{|_{B(0,1]}}=id_{B(0,1]} $, then $\varphi=id_{\U}$.
\end{thm}

\noindent {\bf Proof.} Say that $A \subseteq \U$ is a \textit{set
of uniqueness} iff
$$\forall x,y \in \U \bigg ( \big( \forall z \in A \ d(x,z)=d(y,z) \big) \Rightarrow x=y \bigg ).$$
To prove theorem \ref{regul}, we only need to prove that nonempty
balls of $\U$ are sets of uniqueness: indeed, admit this for a
moment and suppose that $\varphi\colon \U \to \U$ is an isometric
map such that $\varphi_{|_{B(0,1]}}=id_{B(0,1]} $. \\
Let then $x \in \U$: we have
$d(x,z)=d(\varphi(x),\varphi(z))=d(\varphi(x),z)$ for all $z\in
B(0,1]$, so that $\varphi(x)=x$, and we are
done. $\hfill \lozenge$\\

Of course, if $A\subset B$ and $A$ is a set of uniqueness, then
$B$ is one too; therefore, the following proposition is more than
what is needed to prove theorem \ref{regul}:

\begin{prop} \label{uniq}
Let $x_1, \ldots, x_n \in \U$; say that  $f \in E(\{x_1, \ldots, x_n\})$ is nice if
$$\forall i \neq j \ |f(x_i)-f(x_j) |<d(x_i,x_j) \mbox{ and } f(x_i)+f(x_j) >d(x_i,x_j) \ .$$
Then, if $f$ is nice,  $K=\{x_1,\ldots x_n\} \cup \{z \in \U
\colon \forall i \ d(z,x_i)=f(x_i)\}$ is a set of uniqueness.
\end{prop}

\noindent {\bf Proof of Proposition \ref{uniq}.} \\
Let $x \neq y \in \U$; we want to prove that there is some $z \in
K$
such that $d(x,z) \neq d(x,y)$. \\
We may of course assume that $d(x,x_i)=d(y,x_i)$ for all $i$. Let
now $g \in E(\{x_1,\ldots, x_n\} \cup\{x\} \cup\{y\})$ be the
Kat\v{e}tov extension of $f$; notice that $g(x)=g(y)$. \\
Now, pick $\alpha>0$ and define a map $g_{\alpha}$ by: \vspace{0.2cm}\\
- $ g_{\alpha}(x_i) =g(x_i)$ for all $i$,\\
- $g_{\alpha}(y)=g(y)$, and $g_{\alpha}(x)=g(x)-\alpha$. \vspace{0.2cm}\\

\noindent Our hypothesis on $f$ ensures that, if $\alpha >0$ is
small enough, then $g_{\alpha} \in E(\{x_1,\ldots, x_n \}
\cup\{x\} \cup\{y\})$.\\
Hence there is some $z \in \U$ which has the prescribed distances
to $x_1,\ldots x_n, x,y$, so that $z \in K$ and $d(z,x) \neq
d(z,y) .\hfill
\lozenge$\\

\noindent {\bf Remark:} Geometrically, this means that if
$S_1,\ldots S_n$ are spheres of center $x_1,\ldots x_n$, no two of
which are tangent (inwardly or outwardly), and $\cap S_i \neq
\emptyset$, then $\cap S_i \cup\{x_1,\ldots,x_n\}$ is a set of
uniqueness. \\
One may also notice that actually any nonempty sphere is a set of uniqueness. \\
Other examples of sets of uniqueness include the sets $Med(a,b)
\cup\{a,b\}$, where $Med(a,b)=\{z \in \U \colon d(z,a)=d(z,b)\}$
(the proof is similar to the one above); in fact $Med(a,b) \cup
\{a\}$ is a set of uniqueness, whereas $Med(a,b)$ obviously is
not!\\
Also, one may wonder whether the condition in the statement of
Proposition \ref{uniq} is necessary; to see that one needs a
condition of that kind, consider the following example: let
$x_0,x_1$ be any two points such that $d(x_0,x_1)=1$, and let $f$
be defined by $f(x_1)=1$, $f(x_2)=2$. Then, for any point $x$ such
that $d(x,x_0)=d(x,x_1)=\frac{1}{2}$, one necessarily has $f(x)=
\frac{3}{2}$, which proves that the result of Proposition
\ref{uniq} is not true in that case.\\

\noindent Theorem \ref{regul} shows that elements of $Iso(\U)$
have some regularity properties; in particular, if an isometric
map $\varphi$ coincides on an open ball with an isometry $\psi$,
then actually $\varphi=\psi$. One might then wonder, if
$\varphi,\psi \colon \U \to \U$ are two isometric maps such that
$\varphi_{|_B}=\psi_{|_B}$ for a nonempty
ball $B$, whether one must have $\varphi=\psi$. \\
It is easy to see that this is the case if $\varphi(B)=\psi(B)$ is
a set of uniqueness; on the other hand, it is not true in general,
which is the content of the next proposition.

\begin{prop}
There are two isometric maps $\varphi,\psi \colon U \to \U$ such
that $\varphi(x)=\psi(x)$ for all $x \in B(0,1]$, and $\varphi(\U)
\cap \psi(\U)=\varphi(B(0,1])=\psi(B(0,1])$.
\end{prop}
\noindent {\bf Proof. } \\
This result is a consequence of the universality of $\U$: let $X$
denote the metric almagam of two copies of $\U$ (say, $X_1$ and
$X_2$) over $B(0,1]$, and let $\varphi_0$ be an isometry of $X=X_1
\cup X_2$ such that $\varphi_0(X_1)=X_2$, $\varphi_0^2=id$ and
$\varphi_0(x)=x$ for all $x \in B(0,1]$. \\
Pick an isometric embedding $\varphi_1 \colon X \to \U$, and let
$y_0= \varphi_1(0)$; also, let $\eta$ be an isometry from $\U$
onto $X_1$, and let $x_0=\eta^{-1}(0)$. \\
Now let $\varphi=\varphi_1 \circ \eta$, and $\psi=\varphi_1 \circ
\varphi_0 \circ \eta$; by definition of $\varphi_0$, $\varphi$ and
$\psi$ are equal on
$\eta^{-1}(B(0,1])=B(x_0,1]$. \\
Also, one has that \\
$\varphi(\U)=\varphi_1(X_1)$ and $g(\U)=\varphi_1(X_2)$, so
$\varphi(\U) \cap \psi(\U)=\varphi_1(X_1 \cap X_2)= \varphi_1(B(0,1])=
\varphi(B(x_0,1])=\psi(B(x_0,1])\ . \hfill \lozenge $\\

\noindent In a way, the preceding proposition illustrates the fact
that $\U$ contains many non-trivial isometric copies of itself
(other examples include the sets $Med(x_1,\ldots x_n)=\{z \in \U
\colon \forall i,j \ d(z,x_i)=d(z,x_j)\}$). \\
Still, all the isometric copies of $\U$ which we have seen so far
are of empty interior. The next theorem (the proof of which is
based on an idea of Pestov) shows that this is not always the
case:

\begin{thm} \label{nonregular}
If $X \subseteq \U$ is closed and Heine-Borel (with the induced
metric), $M>0$, then $\{z \in \U \colon d(z,X)\geq M\}$ is
isometric to $\U$.
\end{thm}
\noindent (Recall that a Polish metric space $X$ is \textit{Heine-Borel} iff closed bounded balls in $X$ are compact).\\
\noindent In particular, $\U$ and $\U \setminus B(0,1[$ are isometric.\\

\noindent {\bf Proof of Theorem \ref{nonregular}. }\\
We will first prove the result supposing
that $X$ is compact.\\
Define then $Y=\{z \in \U \colon d(z,X)\geq M\}$; $Y$ is a closed
subset of $\U$, so to show that it is isometric to $\U$ we only
need to prove that $Y$ is finitely
injective.\\
Let $y_1,\ldots y_n \in Y$ and $f\in E(\{y_1,\ldots,y_m\})$. Then
there exists a point $c \in \U$ such that $d(c,y_i)=f(y_i)$ for
all $i$; the problem is that we cannot be sure a priori that
$d(c,X) \geq M$.
\\
To achieve this, first define $\ds{\varepsilon=
\min\{f(y_i) \colon 1 \leq i \leq p\} }$. \\
We may of course assume $\varepsilon >0$.\\
$X$ is compact, so we may find $x_1,\ldots x_p \in X$ such that \\
$$\forall x \in X \, \exists j \,\leq p \ \ d(x,x_j) \leq
\varepsilon$$ Let then $g$ be the Kat\v{e}tov extension of $f$ to
$\{y_1,\ldots y_n\} \cup \{x_1,\ldots x_p\}$.
\vspace{0.2cm}\\
By the finite injectivity of $\U$, there is $c \in \U$ such that
$d(c,y_i)=g(y_i)$ for all $i \leq n$ and
$d(c,x_j)=g(x_j)=d(x_j,y_{i_j})+f(y_{i_j}) \geq M+
\varepsilon$ for all $j \leq p$.\\
Since for all $x\in X$, there is $j\leq p$ such that $d(x,x_j)
\leq \varepsilon$, the triangle inequality shows that $d(c,x) \geq
d(c,x_j)-d(x_j,x) \geq M$, hence $c \in Y$, which proves that $Y$
is finitely injective.\vspace{0.2cm}\\
Suppose now that $X$ is Heine-Borel but not compact, and let\\
$Y=\{z \in \U \colon d(z,X)\geq M\}$.
\\
As before, we only need to show that $Y$ is finitely injective; to
that end, let $y_1,\ldots y_n \in Y$ and $f\in
E(\{y_1,\ldots,y_n\})$.\\
Let also $x \in X$ and $m=f(y_1)+d(y_1,x)$. \\
Since $B(x,M+m] \cap X$ is compact, there exists $c \in \U$ such
that $d(c,y_i)=f(y_i)$ for all $i\leq n$, and $d(c,B(x,M+m)) \geq
M$.\\
Then we claim that for all $x' \in X$ we have $d(c,x') \geq M$: if
$d(x',x) \leq M+m$ then this is true by definition of $c$, and if
$d(x',x) >M+m$ then one has $d(c,x')\geq d(x,x')-d(c,x) $, so that
$d(c,x') > M$ (since $d(c,x) \leq f(y_1)+d(y_1,x)=m$). $\hfill
\lozenge$\\

\noindent From the combination of theorems \ref{regul} and \ref{nonregular}, one can easily deduce that:
\begin{coro}\label{debase}
There is an isometry $\varphi$ of $B(0,1]$ such that no isometry
of $\U$ coincides with $\varphi$ on $B(0,1]$.
\end{coro}
\noindent To derive corollary \ref{debase} from the previous
results, let $\varphi \colon \U \to \U \setminus B(0,1[$ be an
isometry, and choose $x \not \in B(0,2]$. There exists, because of
the homogeneity of $\U \setminus B(0,1[$, an isometry $\psi$ of
$\U \setminus B(0,1[$ such that $\psi(\varphi(x))=x$. Thus,
composing if necessary $\varphi$ with $\psi$, we may suppose that
$x$ is a fixed point of $\varphi$. But then $\varphi$ must send
the ball of center $x$ and radius $1$ (in $\U$) onto the ball of
center $x$ and radius $1$ (in $\U \setminus B(0,1[$).\\
Since by choice of $x$ both balls are the same, we see that
$\varphi_{|_{B(x,1]}}$ is an isometry of $B(x,1]$, yet theorem
\ref{regul} shows that no isometry of $\U$ can coincide with
$\varphi$ on $B(x,1]$.$ \hfill \lozenge$\\

\noindent (Using finite injectivity and automatic continuity of
Baire measurable morphisms between Polish groups, one can give a
direct, if somewhat longer, proof of corollary \ref{debase}).
\end{subsection}

\begin{subsection}{Fixed point of isometries} $ $\\
Here we use the tools introduced above - most notably Kat\v{e}tov
maps and the compact injectivity of $\U$ - in order to study some
properties of fixed points of elements of $Iso(\U)$. If $\varphi
\in Iso(\U)$, we let $Fix(\varphi)$ denote its set of fixed
points. \\
Since the isometry class of $Fix(\varphi)$ is an invariant of the
conjugacy class of $\varphi $, one may hope to glean some
information about the conjugacy relation by the study of fixed
points. \\
Clemens, quoted by Pestov in \cite{pest}, conjectured that this
invariant was the weakest possible: the exact content of his
conjecture was that, if $\varphi \in Iso(\U)$, then the set of
fixed points of $\varphi$ is either empty or isometric to $\U$. \\
This turns out to be false in the general case, as we will see
below; this will enable us to provide a lower bound for
the complexity of the conjugacy relation. \\
First, we prove the rather surprising fact that the conjecture
holds for all isometries of finite order (and even for isometries
with totally bounded orbits); so, studying their fixed points will
tell us nothing about, say,
conjugacy of isometric involutions. \\

\noindent We wish to attract the attention of the reader to a
consequence of the triangle inequality, which, though obvious, is
crucial in the following constructions:
$$\forall z \in \U \ \forall x \in \U \ d(z,\varphi(z)) \leq d(z,x)+d(z,\varphi(x)). $$
If $\varphi \colon \U \to \U$ is an isometry, and $x \in \U$, we
let $\rho_{\varphi}(x)=\mbox{diam }\{\varphi^n(x)\}_{n\in \Z}$;
when there is no risk of confusion we simply write it
$\rho(x)$.\vspace{0.2cm}

\begin{lem}\label{migration}
Let $ x_1,\ldots, x_m \in \U$, $f \in E(\{x_1,\ldots,x_m\})$, and
$z \in \U$. Assume that $min\{f(x_i)\} \geq 2 \rho_{\varphi}(z)>0$.\\
Then define \\
$A=\{1\leq i \leq m \colon d(z,x_i) <  f(x_i)-
\frac{\rho_{\varphi}(z)}{2}\}$, $B=\{1\leq i \leq m \colon
d(z,x_i) >  f(x_i) + \frac{\rho_{\varphi}(z)}{2}\}$, and $C=
\{1\leq i \leq m \colon |d(z,x_i)-f(x_i)| \leq \frac{\rho_{\varphi}(z)}{2}\} $. \\
These equations define a Kat\v{e}tov map on $\{\varphi^n(z)\}_{n \in \Z} \cup \{x_i\}_{1\leq i \leq n}$ :\\
- $\forall n \in \Z \  g(\varphi^n(z))=\frac{\rho_{\varphi}(z)}{2}$,\\
- $\forall i \in A \ g(x_i)=d(z,x_i)+\frac{\rho_{\varphi}(z)}{2} $, \\
- $\forall i \in B \ g(x_i)=d(z,x_i)-\frac{\rho_{\varphi}(z)}{2} $, \\
- $\forall i \in C \ g(x_i)=f(x_i)$. \\
Hence, if the orbit of $z$ is totally bounded, there exists $z'
\in \U$ with the prescribed distances to $\{\varphi^n(z)\}_{n \in
\Z} \cup \{x_i\}_{1\leq i \leq n}$;
notice that $\rho(z')\leq \rho(z)$.\\
\end{lem}

\noindent {\bf Proof of lemma \ref{migration}.}\\
To simplify notation, we let $\rho=\rho_{\varphi}(z)$. To check
that the above equations define a Kat\v{e}tov map,
we begin by checking that $g$ is $1$-Lipschitz:\\
First, we have that $|g(x_i)-g(\varphi^n(z))| = |d(z,x_i) + \alpha
-\rho |$, where $|\alpha| \leq \rho$. If $\alpha = \rho$ there is
nothing to prove, otherwise it means that $d(z,x_i) \geq f(x_i)
-\rho $, so that $d(z,x_i) \geq \rho$, which is enough to show
that
$|d(z,x_i) + \alpha -\rho | \leq d(z,x_i)=d(\varphi^n(z),x_i)$.\\
We now let $1 \leq i,j \leq m$ and assume w.l.o.g that
$|g(x_i)-g(x_j)|=g(x_i)-g(x_j)$; the only non-trivial cases are
the following: \vspace{0.2cm}\\
\noindent (a) $g(x_i) = d(z,x_i)+ \alpha$, $g(x_j) = d(z,x_j)+ \beta$, with $\alpha > \beta \geq 0$.  \\
Then one must have $g(x_j)=f(x_j)$, and also $g(x_i) \leq f(x_i)$,
so that $g(x_i)-g(x_j) \leq f(x_i)-f(x_j) \leq d(x_i,x_j)$. \vspace{0.2cm}\\
\noindent (b) $g(x_i)= d(z,x_i)+ \alpha$, $g(x_j)=d(z,x_j)-\beta$,
$0 \leq \alpha,\beta \leq \rho$. Then the definition of $g$
ensures that $g(x_i) \leq f(x_i)$ and $g(x_j) \geq f(x_j)$, so
that
$g(x_i)-g(x_j) \leq f(x_i)-f(x_j) \leq d(x_i,x_j)$.\vspace{0.2cm}\\
\noindent (c)$g(x_i)= d(z,x_i)- \alpha$, $g(x_j)=d(z,x_j)-\beta$,
$0\leq \alpha < \beta$. \\
Then we have $g(x_i)=f(x_i)$, and $g(x_j) \geq f(x_j)$, so
$g(x_i)-g(x_j) \leq f(x_i)-f(x_j)$. \\

\noindent We proceed to check the remaining inequalities: \\

\noindent - $g(\varphi^n(z))+g(\varphi^m(z))=2 \rho \geq d(\varphi^n(z),\varphi^m(z))$ by definition of $\rho$;\\
- $g(\varphi^n(z))+g(x_i)= \rho + d(z,x_i) + \alpha$, where
$|\alpha| \leq \rho$, so
$g(\varphi^n(z))+g(x_i) \geq d(z,x_i)=d(\varphi^n(z),x_i)$. \\
The last remaining inequalities to examine are that involving
$x_i,x_j$; we again break the proof in subcases, of which only two
are not trivial:  \vspace{0.2cm}\\
\noindent (a) $ g(x_i)=d(z,x_i) + \alpha $ and
$g(x_j)=d(z,x_j)-\beta$, where $0 \leq \alpha < \beta$. Then
$g(x_i)=f(x_i)$, and $g(x_j)
\geq f(x_j)$, so that $g(x_i)+g(x_j) \geq d(x_i,x_j)$. \vspace{0.2cm}\\
\noindent (b) $g(x_i)=d(z,x_i)- \alpha$, $g(x_j)=d(z,x_j)-\beta$:
then we have both that $g(x_i) \geq f(x_i)$ and $g(x_j) \geq
f(x_j)$, so
we are done. $\hfill \lozenge$\\

\noindent This technical lemma enables us to prove the following
result, which is nearly enough to prove that $Fix(\varphi)$ is
finitely injective:

\begin{lem}\label{lem2bis}
Let $\varphi$ be an isometry of $\U$ with totally bounded orbits,
$x_1,\ldots, x_m \in Fix(\varphi)$, $f \in E(\{x_1,\ldots,x_m\})$,
and $\varepsilon >0$. Then one (or both) of the
following assertions is true: \\
- There exists $z \in \U$ such that $\rho_{\varphi}(z) \leq
\varepsilon$ and $d(z,x_i)=f(x_i)$ for all $i$ \\
- There is $z \in Fix(\varphi)$ such that $|f(x_i)-d(z,x_i)| \leq
\varepsilon$.\\
\end{lem}

\noindent {\bf Proof of lemma \ref{lem2bis}: }\\
Let $x_1,\ldots, x_m \in Fix(\varphi)$, $f \in
E(\{x_1,\ldots,x_m\})$, and $\varepsilon >0$, which  we assume
w.l.o.g to be  strictly smaller than $\min\{f(x_i) \colon
i=1 \ldots n\} $\\
We may assume that
$$\gamma = \inf \big \{\sum_{i=1}^m |f(x_i)-d(x,x_i)| \colon x \in Fix(\varphi) \big \} >0 \ .$$ \\
Let $x \in Fix(\varphi)$ be such that $\sum_{i=1}^m
|f(x_i)-d(x,x_i)| \leq \gamma + \frac{\varepsilon}{4}$. \\
We let $z$ be any point such that\\
- $d(z,x)=\frac{\varepsilon}{2}$; \\
- $\forall i=1,\ldots, m \
|d(x,x_i)-f(x_i)| \leq \frac{\varepsilon}{2}
\Rightarrow d(z,x_i)=f(x_i)$ ;\\
- $\forall i=1,\ldots, m \ d(x,x_i) \geq f(x_i) +
\frac{\varepsilon}{2} \Rightarrow d(z,x_i)=f(x_i) -
\frac{\varepsilon}{2}$ ;\\
- $\forall i=1,\ldots, m \ d(x,x_i) \leq f(x_i) +
\frac{\varepsilon}{2} \Rightarrow d(z,x_i)=f(x_i) +
\frac{\varepsilon}{2}$. \\

\noindent(One checks as above that these equations indeed define a
Kat\v{e}tov map; $z$ cannot be a fixed point of $\varphi$ since it would contradict the definition of $\gamma$,
or the fact that $\gamma >0$)\\
We use lemma \ref{migration} to build a sequence $(z_n)$ of points
of $\U$ such
that: \\
(0) $z_0=z$;\\
(1) $ 0< \rho(z_n) \leq \varepsilon$;\\
(2) $\forall p \in \Z  d(z_{n+1},\varphi^p(z_n)) = \frac{\rho(z_n)}{2}$; \\
(3) $ \forall i \in A_n \ d(z_{n+1},x_i) = d(z_n,x_i) +
\frac{\rho(z_n)}{2}$; \\
(4) $ \forall i \in B_n \ d(z_{n+1},x_i) = d(z_n,x_i) -
\frac{\rho(z_n)}{2} $;\\
(5) $ \forall i \in C_n \ d(z_{n+1},x_i)=f(x_i)$.\\

\noindent Suppose the sequence has been constructed up to rank
$n$: since $\{x_1,\ldots x_m\}, z_n, f$ satisfy the hypothesis of
lemma \ref{migration}, we may find a point $z'$ with the
prescribed distances to $\{\varphi^p(z_n)\} \cup \{x_1,\ldots
x_m\}$. As before, $z'$ cannot be fixed, since it would contradict
the definition of $\gamma$; we let $z_{n+1}=z'$, and the other conditions are all ensured by lemma \ref{migration}.\\

\noindent If we do not obtain in finite time a $z_n$ such that
$\rho(z_n) \leq \varepsilon$ and $d(z_n,x_i)=f(x_i)$ for all $i$,
then either $A_n$ or $B_n$ is nonempty for all $n$; hence (3) and
(4) imply that $\sum
\rho(z_n)$ converges. Therefore, $z_n$ converges to some fixed point $z^{\infty}$. \\
Necessarily, there was some $i$ such that $|d(z_0,x_i)-f(x_i)|
\leq |d(x,x_i)-f(x_i)|-\frac{\varepsilon}{2}$, so $\sum_{i=1}^m
|f(x_i)-d(z_0,x_i)|\leq \gamma -\frac{\varepsilon}{4}$. \\
By construction, $\sum_{i=1}^m |f(x_i)-d(z^{\infty},x_i)| \leq
\sum_{i=1}^m |f(x_i)-d(z_0,x_i)| $, which contradicts the
definition of $\gamma$. $\hfill \lozenge$\\
$ $\\

\noindent This is not quite enough to produce fixed points with
prescribed distances to some finite set of fixed points; the
following lemma ensures that it is indeed possible:

\begin{lem} \label{lem3}
Let $\varphi$ be an isometry of $\U$ with totally bounded orbits,
$x \in \U$ be such that $\rho_{\varphi}(x)\leq 2\varepsilon $
, and assume that $Fix(\varphi) \neq \emptyset$. \\
Then there exists $y \in \U$ such that : \\
- $\forall n \in \Z \ d(y,\varphi^{n}(x))=d(y,x) \leq \varepsilon$ \\
- $\rho_{\varphi}(y) \leq \varepsilon$.
\end{lem}

\noindent {\bf Proof of lemma \ref{lem3}.} \\
Let $x, \varphi$ be as above; let also
$$E=\{y \in \U \colon \forall n \in \Z \ d(y,\varphi^{n}(x))=d(y,x)
\mbox{ and } \rho(y) \leq \varepsilon\}$$

\noindent Notice that $E$ is nonempty, since any fixed point of
$\varphi$ belongs to $E$. \\

\noindent Now let $\alpha=\inf\{d(y,x) \colon y \in E\}$; we want
to prove that $\alpha \leq \varepsilon$. If not, let $\delta >0$
and pick $y
\in E$ such that $d(y,x)<\alpha + \delta$. \\
Let now $\rho(y)=\beta \leq \varepsilon$; one checks as above that
the following map $g$ belongs to
$E(\{\varphi^n(x)\}\cup\{y\})$: \\
- $\forall n \in \Z \ g(\varphi^n(x))=\max(\varepsilon,d(y,x)- \frac{\beta}{2})$. \\
- $\forall n \in \Z \ g(\varphi^n(y))=\frac{\beta}{2}$.\\
Since the orbits of $\varphi$ are totally bounded, there exists $z
\in \U$ with the prescribed distances; consequently $z \in E$,
and we see that necessarily $\beta < 2 \delta$.\\
Letting $\delta $ go to $0$, there are only two cases to consider:\\

\noindent (1) one may find $y \in \U$ such that $ \ d(y,x)=
\alpha$ and $0< \rho(y)) \leq \varepsilon$.
\\
As before, we may find $z$ such that \\
- $\forall n \in \Z \ d(z,\varphi^n(x))=\max(\varepsilon,d(y,x)- \frac{\rho(y)}{2})$. \\
- $\forall n \in \Z \ d(z,\varphi^n(y))=\frac{\rho(y)}{2}$.\\
Notice that $z \in E$, and $ d(z,x) < \alpha$, which is absurd. \vspace{0.2cm}\\
(2) For all $p \in \N^*$ there is a fixed point $y_p$ such that $
\alpha \leq d(y_p,x)<\alpha+\frac{1}{p}$. If so, let $p$ be big
enough that
$\frac{1}{p}< \frac{\varepsilon}{2}$, and consider the following map: \\
- $g(y_p)=\frac{1}{p}$ \\
- $\forall n \in \Z \ g(\varphi^n(x))=d(y_p,x)-\frac{1}{p}$.\\
A direct verification shows that $g \in
E(\{\varphi^n(x)\}\cup\{y_p\}$, therefore there is $z \in \U$ with
the desired distances; to conclude, notice again that $z \in E$
and $ d(z,x) < \alpha .\hfill \lozenge$\\

We have finally done enough to obtain the following result:
\begin{thm}\label{ordrefini}
If $\varphi \colon \U \to \U$ is an isometry whose orbits are
totally bounded, and $Fix(\varphi)$ is nonempty, then
$Fix(\varphi)$ is isometric to $\U$.
\end{thm}

\noindent {\bf Proof.} Recall that a nonempty metric space $X$ is
said to have the \textit{approximate extension property} iff
$$\forall x_1,\ldots, x_n \in X \ \forall f \in E(\{x_1,\ldots, x_n\}) \ \forall \varepsilon >0 \  \exists z \in X \
|d(z,x)-f(x)|\leq \varepsilon\ .$$

\noindent It is a classical result that, up to isometry, $\U$ is
the only complete, nonempty, separable metric space with the
approximate extension property. So, to prove Theorem
\ref{ordrefini}, it is enough to prove that $Fix(\varphi)$ has the
approximate extension
property.\\
To prove this, notice first that lemma \ref{lem3} implies that,
for all $x \in X$ such that $\rho_{\varphi}(x) \leq \varepsilon$,
there is a fixed point $y$ such that $d(y,x) \leq 2 \varepsilon$. \\
Let now $x_1,\ldots, x_n \in Fix(\varphi)$,$f \in E(\{x_1,\ldots,
x_n\})$, and $\varepsilon >0$. \\
Lemma \ref{lem2bis} tells us that : \\
- there exists a point $z$ such that $\rho_{\varphi}(z) \leq
\frac{\varepsilon}{2}$, and
$d(z,x_i) =f(x_i)$ for all $i=1 \ldots n$, or \\
- there exists $z \in Fix(\varphi)$ such that
$|d(z,x_i)-f(x_i)|\leq \varepsilon$ for all $1 \leq i \leq m$.
\\
In the second case, we have what we wanted; so suppose we are
dealing with the first case, and pick any fixed point $y$ such
that $d(y,z) \leq \varepsilon$. Then $y \in Fix(\varphi)$, and
$|d(y,x_i)-f(x_i)|\leq
\varepsilon$ for all $1\leq i \leq n \ . \hfill \lozenge$ \\

\noindent It turns out that the situation is very different when
it comes to studying isometries with non totally bounded orbits;
one may still prove, using the same methods as above, that if
$\varphi$ is an isometry with a fixed point $x$, then on any
sphere $S$ centered in $x$ and for any $\varepsilon >0$ there is
$z \in S$ such that $d(z,\varphi(z)) \leq \varepsilon$. This is
not enough to ensure the existence of other fixed points than $x$.

\begin{thm}
Let $X$ be a Polish metric space.\\
There exists an isometric copy $X' \subset \U$ of $X$, and an
isometry $\varphi$ of $\U$, such that $Fix(\varphi)=X'$.
\end{thm}

\noindent {\bf Proof .} \\
We may of course assume that $X \neq \emptyset$.\\
We first need a few definitions: if $X$ is a metric space, we
denote by $E(X,\omega,\Q)$ the set of functions $f \in
E(X,\omega)$ which take rational values on their support (This set
is countable if $X$ is). \\
Also, if $X_0 \subset X$ are two countable metric spaces, and
$\varphi$ is an isometry of $X$, we want to find a condition on
$(X,X_0,\varphi)$ which expresses the idea that \\
"$\varphi$ fixes all the points of $X_0$, and for each $x \in X
\setminus X_0$, $\varphi^n(x)$ gets to be as far away from $x$ as
possible". The following definition is a possible way to translate
this naive idea into formal mathematical language:\\

\noindent We say that $(X,X_0,\varphi)$ has property (*) if: \\
- $\forall x \in X_0 \ \varphi(x)=x$.\\
- $\forall x_1,x_2 \in X\ \liminf_{|p| \to + \infty}
d(x_1,\varphi^p(x_2)) \geq d(x_1,X_0)+d(x_2,X_0)$.\\

\noindent The following lemma, which shows that this property is
suitable for an inductive construction similar to Kat\v{e}tov's,
is the core of the proof:

\begin{lem}\label{fratfrat} Let $(X,X_0, \varphi)$ have property (*). \\
Then there exists a countable metric space $X'$ and an isometry
$\varphi'$ of $X'$ such that :\\
- $X$ embeds in $X'$, and $\varphi'$ extends $\varphi$.\\
- $\forall f \in E(X,\omega,\Q)\ \exists x' \in X' \ \forall x \in X \ d(x',x)=f(x)$.\\
- $(X',X_0,\varphi_f)$ has property (*) (identifying $X_0$ to its
image via the isometric embedding of $X$ in $X'$).
\end{lem}

\noindent Admit this lemma for a moment; now, let $X_0$ be any
dense countable subset of $X$, and $\varphi_0=id_{X_0}$. Then
$(X_0,X_0,\varphi_0)$ has property (*), so lemma \ref{fratfrat}
shows that we may define inductively countable metric spaces $X_i$
and isometries $\varphi_i \colon X_i \to X_i$
such that:\\

\noindent -$X_i$ embeds isometrically in $X_{i+1}$, $\varphi_{i+1}$ extends $\varphi_i$;\ \\
-$(X_i,X_0,\varphi_i)$ has property (*);\\
-$\forall f \in E(X_i,\omega,\Q) \ \exists z \in X_{i+1}\ \forall
x \in X_i \ d(z,x)=f(x)$. \\

\noindent Let $Y$ denote the completion of $\cup X_i$, and
$\varphi$ be the extension to $Y$ of the map defined by
$\varphi(x)=\varphi_i(x)$ for all $x \in X_i$. \\
By construction, $Y$ has the approximate extension property; since
$Y$ is complete and nonempty, this shows that $Y$ is isometric to
$\U$. \\
The construction also ensures that all points of $X_0$ are fixed
points of $\varphi$, and $ \liminf_{|p| \to + \infty}
d(y_1,\varphi^p(y_2))\geq d(y_1,X_0)+d(y_2,X_0)$ for all $y_1,y_2 \in Y$. \\
Therefore, $Fix(\varphi)$ is the closure of $X_0$ in $\U$; hence
it is isometric to the completion of $X_0$,
so it is isometric to $X$. $\hfill \lozenge$\\

\noindent {\bf Proof of Lemma \ref{fratfrat}.} \vspace{0.2cm}\\
First, let $f \in E(X,\omega, \Q)$;  we let $X(f)=X \cup
\{y^f_i\}_{i \in \Z}$ and define a distance on
$X(f)$, which extends the distance on $X$, by: \vspace{0.2cm}\\
-$d(x,y^f_i)=f(\varphi^{-i}(x))$; \\
-$d(y^f_i,y^f_j)= \inf_{x \in X} (d(y^f_i,x)+d(y^f_j,x))$.\vspace{0.2cm}\\
(In other words, $X(f)$ is the metric amalgam of the spaces $X
\cup\{f \circ \varphi^i\}$ over $X$. )\\
Let $\varphi_f$ be defined by $\varphi_f(y^f_i)=y^f_{i+1}$,
$\varphi_f(x)=\varphi(x)$ for $x\in X$. \vspace{0.2cm}\\
Notice that, by definition of $d$, $\varphi_f$ is an
isometry of $X(f)$, which extends $\varphi$.\\
We claim  that $(X(f),X_0,\varphi_f)$ has property (*). \\
To prove this, let $y,y' \in X(f)$; we want to prove that
$$\liminf_{|p| \to + \infty} d((\varphi_f)^p(y'),y) \geq d(y,X_0)+d(y',X_0)\ .$$
If both $y$ and $y'$ are in $X$, there is nothing to prove. Two
cases remain: \vspace{0.2cm}\\
\noindent (1) $y\in X$, $y'= y^f_j$. Without loss of generality,
we may assume that $j=0$. By definition, we know that
$$d((\varphi_f)^p(y^f_0),y)=f(\varphi^{-p}(y))= \min_{i=1 \ldots n}
(f(x_i)+d(y,\varphi^p(x_i)))$$ for some $x_1,\ldots x_n \in X$ (recall that $f \in E(X,\omega,\Q)$).\\
Let $\varepsilon >0$;  for $|p|$ big enough,
$d(y,\varphi^p(x_i)) \geq d(y,X_0)+d(x_i,X_0) - \varepsilon$.\\
We then have $d((\varphi_f)^p(y^f_0),y) \geq \min_{i=1 \ldots n}
(f(x_i)+d(y,X_0)+d(x_i,X_0) - \varepsilon)$, so
$d((\varphi_f)^p(y^f_0),y) \geq d(y,X_0) + \min_{i=1 \ldots
n}(f(x_i) +d(x_i,X_0)) -
\varepsilon $. \\
Hence $ d((\varphi_f)^p(y^f_0),y) \geq d(y,X_0) +d(y^f_0,X_0)-
\varepsilon ,$ and we are done. \vspace{0.2cm} \\
\noindent (2) $y=y^f_i$ and $y'=y^f_j$; we may assume that $i=0$. \\
Then we have $d(\varphi_f^p(y'),y)= \inf_{x \in X}
(f(x)+f(\varphi^{-p-i}(x))$. \\
We want to prove that
$$\liminf_{|p| \to + \infty} \inf_{x \in X}
(f(x)+f(\varphi^{-p}(x)) \geq 2\inf_{x \in X_0} f(x) \ .
$$
Assume again that $f$ is controlled by $\{x_1,\ldots,x_n\}$,
choose $\varepsilon >0$, and let $|p|$ be big enough that
$d(x_i,\varphi^p(x_j)) \geq
d(x_i,X_0)+d(x_j,X_0)-\varepsilon$ for all $i,j$. \\
Then we have, for all $x \in X$:\\
$f(x)+f(\varphi^{-p}(x))=f(x_i)+d(x,x_i)+f(x_j)+d(x,\varphi^p(x_j))$
for some $i,j$.
\\
Since $d(x,x_i)+d(x,\varphi^p(x_j)) \geq d(x_i,\varphi^p(x_j))$,
we see that there is some $(i,j)$ such that $\inf_{x \in
X}(f(x)+f(\varphi^{-p}(x))=
f(x_i)+d(x_i,\varphi^p(x_j))+f(x_j)$.\\
We know that $d(x_i,\varphi^p(x_j)) \geq
d(x_i,X_0)+d(x_j,X_0)-\varepsilon$, so
$$ \inf_{x \in X}(f(x)+f(\varphi^{-p}(x)) \geq f(x_i)+d(x_i,X_0)+d(x_j,X_0)+f(x_j)-\varepsilon
\geq 2\inf_{x \in X_0} f(x)- \varepsilon .$$

\noindent This is enough to prove that $(X(f),X_0,\varphi_f)$ has
property (*).  \\

\noindent Now, let $X'$ denote the metric amalgam of the spaces
$X(f)$ over $X$, where $f$ varies over $E(X,\omega,\Q)$. It is
countable, and letting $\varphi'(x)= \varphi_f(x)$ for all $X \in
X^f$ defines an isometry of $X'$ which extends $\varphi$. \\
If $f \neq g \in E(X,\omega,\Q)$, let $k_f(g)$ denote the
Kat\v{e}tov of $g$ to $X(f)$; by definition, the metric amalgam of
$X(f)$ and $X(g)$ over $X$ is exactly the space
$(X(f))(k_f(g))$, which is enough to show that $(X',X_0,\varphi)$ has property (*). $  \hfill \lozenge$\\

\noindent This construction has an additional interest, since it
provides a lower bound for the complexity of conjugacy between
isometries of $\U$. Indeed, we may endow any countable graph with
the graph distance, turning it into a countable Polish metric
space; two graphs are isomorphic if, and only if, the
corresponding metric spaces are isometric.\\
Now, let $X$ and $X'$ denote two isometric countable Polish metric
spaces. Let $X_{\infty}=\cup X_i$ and $X'_{\infty}=\cup X'_i$
denote the spaces obtained by our construction, and
$\varphi_{\infty}, \varphi'_{\infty}$ the corresponding
isometries. It is not hard to see that the isometry between $X$
and $X'$ extends to an isometry $\psi \colon X_{\infty} \to
X'_{\infty}$ such that $\psi \circ
\varphi_{\infty}= \varphi'_{\infty} \circ \psi$. \\
Since the completions of $X_{\infty}$ and $X'_{\infty}$ are both
isometric to $\U$, this means that we may, reasoning as in
\cite{gaokec}, build a Borel map \\$ \Psi \colon GRAPH \to
Iso(\U)$ such that any graph $G$ is isometric to $Fix(\Psi(G))$,
and $\Psi(G)$ and $\Psi(G')$ are conjugate if $G$ and $G'$ are
isomorphic. \\
Conversely, assume that there is $\varphi \in Iso(\U)$ such that
$\varphi \circ \Psi(G)= \Psi(G') \circ \varphi$; this implies that
$\varphi(Fix(\Psi(G)))=Fix(\Psi(G'))$, and this proves that $G$
and $G'$ are isometric. We have just proved the following result:

\begin{thm}
Graph isomorphism Borel reduces to conjugacy of isometries in
$\U$.
\end{thm}




\end{subsection}

\end{section}


\begin{section}{Trying to extend finite homogeneity} \label{trans}
\begin{subsection}{Reformulating the problem}$ $\\
\noindent The remainder of this article will be devoted to proving
the following result:
\begin{thm} \label{final} Let $X$ be a Polish metric space. The following
assertions are equivalent: \\
(a) X is compact.\\
(b) If $ X_1, X_2 \subseteq \U$  are both isometric to $X$ and $\varphi \colon X_1 \to X_2$ is an isometry,
then there exists $ \tilde{ \varphi} \in Iso(\U)$ which extends $\varphi$.\\
(c) If $ X_1, X_2 \subseteq \U$  are both isometric to $X$, then
there exists $\varphi \in Iso(\U)$ such that $\varphi(X_1)=X_2$.\\
(d) If $X_1 \subseteq \U$ is isometric to $X$ and $f \in E(X_1)$,
there exists $z \in \U$ such that $d(z,x)=f(x)$ for all $x \in
X_1$.\\

\end{thm}

\noindent  $(a) \Rightarrow (b)$ is well-known, as explained in
the introduction (see \cite{huhu} for a proof); $(b) \Rightarrow (c)$ is trivial.\\
To see that $(c) \Rightarrow (d)$, let $X$ have property (c) and
be embedded in $\U$, and $f \in  E(X)$;  the metric space $X_f=X
\cup \{f\}$ embeds in $\U$, so that there exists an isometric copy
$Y=X' \cup\{z\}
\subset \U$ of $X_f$, where $X'$ is an isometric copy of $X$. \\
Notice that, since there exists a copy of $X$ which is
$g$-embedded in $\U$, and all isometric copies of $X$ are
isometric by an isometry of the whole space, all the isometric
copies of $X$ are necessarily
$g$-embedded in $\U$. \\
Let now $\varphi$ be the isometry from $X$ to $X'$ which sends any
point $x$ to its copy in $X'$; we have, by definition,
$d(z,\varphi(x))=f(x)$. Pick now an isometry $\psi$ of $\U$ which
maps $X' \to X$; then $d(\psi(z),\psi \circ \varphi(x))=f(x))$.\\
Consequently, if we let $\rho$ be an isometry of $\U$ which
extends $(\psi \circ \varphi) ^{-1}$ then we have for all $x \in
X$:
$$d(\rho(\psi(z)),x)=d(\psi(z),\rho^{-1}(x))=d(\psi(z),\psi \circ \varphi(x))=d(z,\varphi(x))=f(x). $$
So $\rho(\psi(z))=z'$ is such that $d(z',x)=f(x)$ for all $x \in
X$, and $X$ has property $(d)$.\\

\noindent It only remains to show that $(c) \Rightarrow (a)$;
this turns out to be the hard part of the proof. \\
\noindent If $X \subseteq \U$ is closed, define $\Phi^X \colon \U
\to E(X)$ by
 $\Phi^X(z)(x)=d(z,x)$.\\
Notice that $\Phi^X$ is $1$-Lipschtitzian. Property $(d)$ in
theorem \ref{final} is equivalent to $\Phi^{X_1}$ being onto for
any isometric copy $X_1 \subseteq \U$ of $X$; but $\Phi^{X_1}(\U)$
is necessarily separable since $\U$
is, so we see that for $X$ to have property $(d)$ it is necessary that $E(X)$ be separable. \\

\noindent The next logical step is to determine the Polish metric
spaces $X$ such that $E(X)$ is separable.
\end{subsection}

\begin{subsection}{Compactly tentacular spaces}$ $\\
One can rather easily narrow the study:

\begin{prop}
If $X$ is Polish and not Heine-Borel, then $E(X)$ is not separable.
\end{prop}
\noindent {\bf Proof:} The hypothesis tells us that there exists
$M,\varepsilon >0$ and $(x_i)_{i \in \N}$ such that
$$\forall i \neq j \ \varepsilon \leq d(x_i,x_j) \leq M.$$
If $A \subseteq \N$, define $f_A:\{x_i\}_{i \geq 0} \to \R$ by
$f_A(x_i)= \begin{cases}
M &\text{if}\;\;i \not \in A\\
M+\varepsilon &\text{else}
\end{cases} $  . \\
It is easy to check that for all $A \subseteq \N$, $f_A \in
E(\{x_i\}_{i \geq 0})$,
and if $A \neq B$ one has $d(f_A,f_B)= \varepsilon$ (where $d$ is the distance on $E(\{x_i\}$) ). \\
Hence $E(\{x_i\}_{i \geq 0})$ is not separable; since it is
isometric to a subspace of $E(X)$ (see section \ref{nota}), this
concludes the proof. $\hfill
\lozenge$\\

\noindent So we know now that, to have property $(d)$ of theorem
\ref{final}, a metric space $X$ has to be Heine-Borel; at this
point, one could hope that either only compact sets are such that
$E(X)$ is separable, or all Heine-Borel Polish spaces have this
property. Unfortunately, the situation is not quite so simple, as
the following two examples show:
\begin{example} If $\N $ is endowed with its usual distance, then
$E(\N)=\overline{E(\N,\omega)}$.
\end{example}

\noindent Indeed, let $f \in E(\N)$; then one has for all $n$ that
$|f(n)-n| \leq f(0)$, and also $f(n+1) \leq
f(n)+1$. This last inequality can be rewritten as $f(n+1)-(n+1)\leq f(n)-n$.\\
So $f(n)-n$ converges to some $a \in \R$; let $\varepsilon >0$ and
choose $M$ big enough that $n \geq M \Rightarrow |f(n)-n-a|\leq \varepsilon$. \\
Then, for all $n \geq M$, one has
$$0 \leq f(M)+n-M-f(n)=(f(M)-M-a)-(f(n)-n-a) \leq 2\varepsilon .$$
If one lets, for all $i$, $f_i$ be the Kat\v{e}tov extension of
$f_{|_{[0,i]}}$, then \\ $f_i \in E(\N,\omega)$ and we have just
shown that $(f_i)$ converges uniformly to $f$.
\\
Replacing the sequence $f(n)-n$ by the function $f(x)-x$, one
would have obtained the same result for any subset of $\R$
(endowed with its usual metric, of course); actually, one may use
the same method to prove that $E(\R^n,||.||_{1})$ is
separable for all $n$. \\

\noindent The situation turns out to be very different when $\R^n$
is endowed with other norms, as the following example shows.

\begin{example}\label{ex2} If $n \geq 2$ and $\R^n$ is endowed with the euclidian distance,
then $E(\R^n)$ is not separable.
\end{example}

\noindent We only need to prove this for $n=2$, since
$E(\R^2,||.||_{2})$ is isometric to a closed subset
of $E(\R^n,||.||_{2})$ for any $n \geq 2$ .\\
Remark first that it is easy to build a sequence $(x_i)$ of points
in $ \R^2$ such that \\
$d(x_{i+1},0)\geq d(x_i,0)+1$ for all $i$, and
$$\forall i>j \in \N ,\ d(x_i,0) \leq d(x_i,x_j)+d(x_j,0)-1 \qquad (*)$$
One can assume that $d(x_i,0) \geq 1$ for all $i$; now define $f
\colon \{x_i\}_{i \geq 0}  \to \R$ by $f(x_i)=d(x_i,0)$.
Obviously, $f$ is a Kat\v{e}tov map. \\
If $A \subseteq \N$ is nonempty, define $f_A \colon \{x_i\}_{i
\geq 0} \to \R$ as the  Kat\v{e}tov extension of $f_{|_{
\{x_i\colon i \in A\}}}$.\\
Suppose now that $A \neq B$ are subsets of $\N$, let $m$ be the
smallest element of
$A \Delta B$, and assume without loss of generality that $m\in A$.\\
Then one has $f_A(x_m)=d(x_m,0)$, and $f_B(x_m)=d(x_m,x_i)+d(x_i,0)$ for some $i \neq
m$.\\
If $i < m$, then $(*)$ shows that $f_B(x_m)-f_A(x_m) \geq 1$; if $i>m$,
then $f_B(x_m)-f_A(x_m) \geq d(x_i,0)-d(x_m,0) \geq 1$.\\
In any case, one obtains $d(f_A,f_B) \geq 1$ for any $A \neq B$, which shows that $E(\{x_i\}_{i \geq 0})$ is not separable.\\
Hence $E(\R^2,|| . ||_2)$ cannot be separable either.\\

\noindent These two examples have something in common: in the
first case, the fact that all points lie on a line gives us that
$\overline{E(X,\omega)}=E(X)$; in the second case, the existence
of an infinite sequence of points on which the triangle inequality
is always far from being an equality enables
us to prove that $E(X)$ is not separable.\\
It turns out that this is a general situation, and we can now
characterize the spaces $X$ such that $E(X)$ is separable:\\
Let $(X,d)$ be a nonempty metric space. \\
For $\varepsilon>0$, a sequence $(u_n)_{n\in\N}$ in $X$ is said to
be {\em $\varepsilon$-good-inline} if for every $r \geq 0$ we have
$ \sum_{i=0}^r d(u_i,u_{i+1})\leq d(u_0,u_{r+1})+\varepsilon$. \\
A sequence $(u_n)_{n \in\N}$ in $X$ is said to be {\em inline} if
for every $\varepsilon>0$ there exists $N\geq 0$ such that
$(u_{n+N})_{n\in\N}$ is $\varepsilon$-good-inline.

\begin{thm}\label{carac} Let $X$ be a Polish metric space.\\
The following assertions are equivalent:\\
(a)$E(X)=\overline{E(X,\omega)}$\\
(b)$E(X)$ is separable\\
(c)$X$ is Heine-Borel and \\
$\qquad \forall \delta >0 \  \forall (x_n) \ \exists N \in \N \
\forall n \geq N \  \exists i \leq N \ d(x_0,x_n)\geq
   d(x_0,x_i)+d(x_i,x_n)-\delta$.\\
(d) Any sequence of points of $X$ admits an inline subsequence.\\

\end{thm}

\noindent{\bf Proof of Theorem \ref{carac}.} \vspace{0.2cm}\\
$(a) \Rightarrow (b)$ is obvious; the proof of $\neg (c)
\Rightarrow \neg (b)$ is similar to Example \ref{ex2},
so we leave it as an exercise for the interested reader.  \\
To see that $(c) \Rightarrow (d)$, one simply needs to repeatedly
apply the pigeon-hole principle.
It remains to prove that $(d) \Rightarrow (a)$.\\
For that, suppose by contradiction that some Polish metric space
$X$ is Heine-Borel, has property (d), but not property (a).\\
Choose then $f \in E(X) \setminus \overline{E(X,\omega)}$, and let
$f_n$ be the Kat\v{e}tov extension to $X$
of $f_{|_{B(z,n]}}$ (where $z$ is some point in $X$).\\
Then for all $x \in X$, $n \leq m$ one has $f_n(x) \geq f_m(x) \geq f(x)$; hence the sequence $(d(f_n,f))$ converges to
some $a \geq 0$. \\
Notice that, since closed balls in $X$ are compact, each $f_n$ is in $\overline{E(X,\omega)}$: this proves that
$a>0$, and one has $d(f_n,f) \geq a$ for all $n$.\\
One can then build inductively a sequence $(x_i)_{i \geq 1}$ of
elements of $X$, such that for all $i \geq 1$ $d(x_{i+1},z) \geq
d(x_i,z)+1$ and
$$f(x_i) \leq \min_{j<i} \{f(x_j)+d(x_i,x_j)\} - \frac{3a}{4}$$
Since $|f(x_i)-d(x_i,z)| \leq f(z)$, one can assume, up to some
extraction, that
$(f(x_i)-d(x_i,z))$ converges to some $l \in \R$. \\
Now, let $\delta= \frac{a}{4}$. (d) tells us that we can extract
from the sequence $(x_i)$ a subsequence $x_{\varphi(i)}$ having
the additional property that
$$ \forall 1\leq j \leq i, \qquad  d(z,x_{\varphi(i)}) \geq d(z,x_{\varphi(j)})+d(x_{\varphi(i)},x_{\varphi(j)})-\delta$$
To simplify notation, we again call that subsequence $(x_i)$.\\
Choose then $M \in \N$ such that $n \geq M \Rightarrow
|f(x_n)-d(x_n,z)-l| \leq \frac{\delta}{2}$.\\
For all $n \geq M$, we have \\
$f(x_M)+d(x_M,x_n)-f(x_n)=(f(x_M)-d(x_M,z)-l) - (f(x_n)-d(x_n,z)-l) + (d(x_M,z)-d(x_n,z)+d(x_M,x_n))$, so that\\
$f(x_M)+d(x_M,x_n)-f(x_n) \leq
2\delta=\frac{a}{2}<\frac{3a}{4}$ .\\
This contradicts the definition of the sequence $(x_i)$, and we
are done. $\hfill \lozenge$\\

\noindent For lack of a better word, we will call (for now...) a
non totally bounded metric space $X$ such that
$\overline{E(X,\omega)}=E(X)$ a \textit{compactly tentacular}
metric space. \\

\noindent It is worth pointing out that in the course of the proof
of theorem \ref{carac}, we proved that, if $X$ is compactly
tentacular and $f \in E(X)$, then for any $\varepsilon >0$ there
exists a compact $K \subseteq X$ such that
$d(f,k(f_{|_K})) < \varepsilon$. \\

The following fact is worh stating:

\begin{thm}
$QL=\{F \in \mathcal{F}(\U) \colon F \mbox{ is compactly
tentacular}\}$ is a Borel subset of $\mathcal{F}(\U)$, endowed
with the Effros Borel structure.
\end{thm}
\end{subsection}


\begin{subsection}{End of the proof of theorem \ref{final}}$ $\\
Now we are ready to finish the proof of theorem \ref{final}; we
need to study the case of compactly tentacular spaces. \\
Let $X$ be a compactly tentacular metric space; we wish to build a
copy $X' \subset \U$ of $X$ such that $\Phi^{X'}(\U) \neq
E(X')$.\\
So, it is natural to try to build an isometric copy $X_1 \subset
\U$ of $X$ such that $\Phi^{X_1}(\U)$ is as small as possible.\\
\noindent To do this, we need a definition:\\
If $X$ is a metric space and $\varepsilon
>0$, we say that $f \in E(X)$ is $\varepsilon$-\textit{saturated}
if there exists a compact $K \subset X$ such that, for any $g \in
E(X)$, $g_{|_K}=f_{|_K} \Rightarrow d(f,g) \leq \varepsilon$. For
convenience, we say that
such a compact $K$ \textit{witnesses} the fact that $f$ is $\varepsilon$-saturated. \\
We say that $f$ is \textit{saturated} if it is
$\varepsilon$-saturated for all $\varepsilon>0$; simple examples
of saturated maps are given by maps of the form $z \mapsto
d(x,z)$, where $x \in X$ (since for any $\varepsilon >0$ one can take $K=\{x\}$). \vspace{0.2cm}\\
A more interesting example is the following: let $X=\N$, and $f\in
E(\N)$ be such that $f(0)=f(1)=1/2$. \\Then the triangle
inequality implies that $f(n+2)=n+3/2$ for all $n \in \N$, which
shows that $f$ is saturated. In other words, such a map is
necessarily contained
in $\Phi^{\N}(\U)$ whenever $\N$ is embedded in $\U$.\\

\noindent It is easy to see that if $X$ is a noncompact metric
space there is $f \in E(X)$ which is not saturated. Thus, the
following proposition is enough to finish the proof of Theorem
\ref{final}:
\begin{prop} \label{laclasse}
Let $X$ be a compactly tentacular space. There exists an isometric
copy $X' \subseteq \U$ of $X$ such that $\Phi^{X'}(z)$ is
saturated for all $z \in \U$.
\end{prop}

\noindent We will use in the proof of Proposition \ref{laclasse}
some simple properties of $\epsilon$-saturated maps on compactly
tentacular metric spaces, which we regroup in the following
technical lemma in the hope of making the proof itself clearer:
\begin{lem} \label{saturated}Let $X$ be a compactly tentacular Polish metric space.\vspace{0.2cm}\\
(1) If $\varepsilon >0$ and $f \in E(X)$ is not
$\varepsilon$-saturated, then for any compact $K \subseteq X$
there is $g \in E(X)$ such that $g_{|_K}=f_{|_K}$ and $g(x)>
f(x)-\varepsilon$ for some $x \in X$. \vspace{0.2cm}\\
(2) If $f \in E(X)$ is saturated, then  for any $\varepsilon
>0$ there exists a compact $K \subseteq X$ such that
$$ \exists M \, \forall x \in X \ d(x,K)\geq M \Rightarrow \exists z
\in K \ f(z)+f(x) \leq d(z,x)+\varepsilon.$$ \vspace{0.2cm}
(3) Let $f_n \in E(X)$ be $\varepsilon_n$-saturated maps such that : \\
- For any $n$ there exists a compact $K_n$ which witnesses the
fact that $f_n$ is $2\varepsilon_n$-saturated, and such that $m
\geq n \Rightarrow {f_m}_{|_{K_n}}= {f_n}_{|_{K_n}}$.\\
- $\varepsilon_n \to 0$.\\
- $\cup K_n=X$ \\
Then $f_n$ converges uniformly to a saturated Kat\v{e}tov map
$f$.

\end{lem}

\noindent {\bf Proof of Lemma \ref{saturated}} \vspace{0.2cm}\\
(1)Since $X$ is compactly tentacular, there exists a compact set
$L$ such that $d(k(f_{|_L}),f) \leq \frac{\varepsilon}{2}$; we may
assume that $K \supseteq L$. \\
Since $f$ is not $\varepsilon$-saturated, we know that there is $g
\in E(X)$ such that $g_{|_K}=f_{|_K}$ and $d(g,f) >
\varepsilon$.\\
Thus there exists $x$ such that $|f(x)-g(x)|> \varepsilon$.\\
Yet, by definition of a Kat\v{e}tov extension, we necessarily have
that  $g \leq k(f_{|_K}) \leq k(f_{|_L}) \leq f+
\frac{\varepsilon}{2}$, so that $|f(x)-g(x)|> \varepsilon$ is only
possible if $f(x)-g(x) > \varepsilon$, i.e
$g(x) < f(x)-\varepsilon$. \\

\noindent(2)Let $f$, $\varepsilon >0$ be as above, and $K$ be a
compact witnessing the fact that $f$ is $\frac{\varepsilon}{2}$-saturated.\\
Now, pick any $x$ such that $d(x,K)\geq M=2\max\{f(x) \colon x \in K \} + \varepsilon$.\\
Suppose by contradiction that one has $f(x)+f(z)> d(z,x)+
\varepsilon$ for any $z \in K$, and let $g$ be defined on $K \cup
\{x\}$ by $g_{|_K}=f_{|_K}$ and
$g(x)=f(x)-\varepsilon$.\\
Then for any $z \in K$ we have \\
$|g(x)-g(z)|=|f(x)-f(z)-\varepsilon|=f(x)-f(z)-\varepsilon \leq
f(x)-f(z) \leq d(x,z)$. Also, for any $z \in K$ one has
$g(x)+g(z)=f(x)+f(z)-\varepsilon >d(z,x) $. \\
Finally, it is obvious that $|g(z_1)-g(z_2)|=|f(z_1)-f(z_2)| \leq
d(z_1,z_2) \leq f(z_1)+f(z_2)=g(z_1)+g(z_2)$ for all $z_1,z_2 \in K$. \\
Consequently, the Kat\v{e}tov extension $k (g)$ of $g$ to $X_p$ is
such that $k(g)_{|_K}=f_{|_K}$ and
$d(f,k(g)) \geq \varepsilon$, which contradicts the definition of $K$. \\

\noindent (3)Let $X$, $f_n$, $\varepsilon_n$ and $K_n$ be as in the statement of \ref{saturated}(3).  \\
Then $(f_n)$ obviously converges pointwise to some Kat\v{e}tov map
$f$, and we have to show that $f$ is saturated and the convergence
is actually uniform.\\
To that end, let $\varepsilon >0$ and choose $N$ such that
$2\varepsilon_N \leq \frac{\varepsilon}{2}$.\\
Then we have, for all $n \geq N$, that
${f_n}_{|_{K_N}}={f_N}_{|_{K_N}}$, which by definition of $K_N$
implies that $d(f_n,f_N) \leq  2\varepsilon_N$. But then one gets
$d(f_n,f_m) \leq \varepsilon$ for any $n,m \geq N$, which shows
that $(f_n)$ is Cauchy, which proves that the convergence is uniform.\\
To show that $f$ is saturated, let again $\varepsilon >0$ and find
$n$ such that $2\varepsilon_n \leq \frac{\varepsilon}{2}$ and
$d(f_n,f) \leq \frac{\varepsilon}{2}$.\\
Then any Kat\v{e}tov map $g$ such that
$g_{|_{K_n}}=f_{|_{K_n}}={f_n}_{|_{K_n}}$ has to satisfy \\$d(f,g)
\leq d(f,f_n)+d(f_n,g) \leq \varepsilon \ .\hfill \lozenge$\\

\noindent{\bf Proof of Proposition \ref{laclasse}.}\\
The method of proof we intend to use is classical (cf section
\ref{nota}): we let $X_0=X$, and define inductively metric spaces
$X_i$ such that $X_{i+1} \supseteq X_i$, $\cup X_i$ is finitely
injective, and $\cup X_i$ has the desired property; then its
completion will
be isometric to $\U$ and satisfy the result of the theorem.\\
So, let as promised $X_0=X$, and define
$$X_{i+1}=
\{ f \in E(X_i) \colon f_{|_{X_0}} \mbox { is saturated }\}\ .$$
(This makes sense since, as in section \ref{nota}, we may
assume, using the Kuratowski map, that $X_i \subseteq X_{i+1}$). \\
As usual, we let $Y$ denote the completion of $\cup X_i$, and need
only prove that $Y$ is finitely injective to conclude the proof.\\
For that, it is enough to show that $\cup X_i$ is finitely
injective; take then $\{x_1,\ldots x_n\} \subseteq X_p$ (for some
$p \geq 0$) and $f \in E(\{x_1,\ldots x_n\})$.\\
We only need to find a map $f \in E (X_p)$ which takes the
prescribed values on $x_1,\ldots x_n$ and whose restriction to
$X_0$ is saturated, since this will belong to $X_{p+1}$ and have
the
desired distances to $x_1,\ldots x_n$.\\
To achieve this, we use the following lemma:

\begin{lem} \label{anarchyintheuk}Let  $x_1,\ldots, x_n \in X_p$, $f \in
E(\{x_1,\ldots, x_n\}$. \\Let also $f'\in E(X_p)$ and $\varepsilon
>0$ be such that $f'(x_i)=f(x_i)$ for all $i$, and $f'_{|_{X_0}}$ is not $\varepsilon$-saturated.\\
Then, for any compact $K_0 \subset X_p$, there exists $g \in
E(X_p)$ such that
$$\forall i=1 \ldots n \ g(x_i)=f(x_i),\
g_{|_{K_0}}=f'_{|_{K_0}} \mbox{ and } \exists x \in X_p \setminus
K_0\ g(x)\leq f'(x)-\frac{\varepsilon}{2}. \qquad (*)$$
\end{lem}

\noindent {\bf Proof of lemma \ref{anarchyintheuk}}\\
To simplify notation below,
fix some point $z_0 \in K_0$.\\
Since $f'$ is not $\varepsilon$-saturated, we can find $y_1 \not
\in K_0$ such that $f'(y_1)+f'(z) > d(y,z)+\varepsilon $ for all
$z \in K_0$. Letting $K_1=B(z_0,2d(z,y_1))$ we can apply
the same process and find $y_2$, and so on. \\
It is not hard to see that one can indefinitely continue this
process, and one can thus build a sequence $(y_n)$ such that
$d(y_n,z_0) \to + \infty$, an increasing sequence of compact sets
$(K_i)$ such that $\cup K_i=X$, and
$$\forall i \geq 1 \ \forall z \in K_{i-1} \ f'(y_i)+f'(z) > d(y_i,z)+\varepsilon \
.$$\\
\noindent {\bf Claim: } If one cannot find a map $g$ as in $(*)$,
then there exists $I$ such that
$$ \forall i \geq I \ \exists k_i \ f'(y_i)+f(x_{k_i}) < d(x_{k{_i}},y_i)+\frac{\varepsilon}{2} \ .\qquad (**)$$ \\
\noindent {\bf Proof }: By contradiction, assume that for all $I$
there exists $i \geq I$ such that $f'(y_i)+f'(x_k) \geq
d(x_k,y_i)+\frac{\varepsilon}{2}$ for
all $k=1\ldots n$. \\
Choose $I$ such that $d(y_I,z_0) \geq \max\{f'(z) \colon z \in
K_0\}+ \frac{\varepsilon}{2}$, $f'(y_i) \geq f'(z)$ for all $z \in
K_0$ and $i \geq I$,  $K_I \supseteq B(z_0,2\mbox{diam}(K_0)]$,
then find $i \geq I$ as above. \\
Define a map $g$ on $\{x_k \}_{k=1\ldots n} \cup{K_0} \cup
\{y_i\}$ by $g(y_i)=f'(y_i)- \frac{\varepsilon}{2}$, $g(x)=f'(x)$ elsewhere. \\
By choice of $i$ and since $\ f'(y_i)+f'(z) \geq
d(y,z)+\frac{\varepsilon}{2}$ for all $z \in K_0$, we see that $g$
is Kat\v{e}tov, and that its Kat\v{e}tov extension $k(g)$ to $X_p$
is such that $k(g)(x_i)=f(x_i),\
 k(g)_{|_{K_0}}=f'_{|_{K_0}}$ and $k(g)(y_i)\leq
f'(y_i)-\frac{\varepsilon}{2}$.\\
This concludes the proof of the claim. \\

\noindent Up to some extraction, we may assume that $k_i=k$ for
all $i \geq i$. By definition of $X_p$, we know that the
restriction to $X_0$ of the map $d(x_k,.)$ is saturated, so lemma
\ref{saturated} shows that there exists $J$ such that
$$\forall j >J\, \exists z \in K_J \ \  d(x_k,z)+d(x_k,y_j) \leq d(z,y_j)+\frac{\varepsilon}{4} \ .$$
Combining this with $(**)$, we obtain, for $j > \max(I,J)$, that
there exists $z \in K_{J} \subseteq K_{j-1}$ such that
$f'(y_j)+f(x_k)+d(x_k,z) \leq d(z,y_j)+\frac{\varepsilon}{2}+ \frac{\varepsilon}{4}$. \\
This in turn implies that $f'(y_j)+f'(z) <d(z,y_j)+\varepsilon$,
which contradicts the definition of the sequence $(y_i). \hfill
\lozenge$ \\

\noindent We are now ready to move on to the last step of the
proof of proposition \ref{laclasse}: \\
First, pick $\{x_1,\ldots x_n\} \subseteq X_p$ (for some $p \geq
0$) and $f \in E(\{x_1,\ldots x_n\})$.\\
We wish to obtain $g \in E(X_p)$ such that
$g(x_i)=f(x_i)$ for all $i$, and $g_{|_{X_0}}$ is saturated.\\
Letting $\varepsilon_0=\inf\{\varepsilon>0 \colon k(f)_{|_{X_0}}
\mbox { is } \varepsilon-\mbox{saturated } \}$, we only need to
deal with the case $\varepsilon_0 >0$ .\\
We have shown that if $k(f)_{|_{X_0}}$ is not
$\varepsilon$-saturated then for any compact $K \subseteq X_p$ we
may find $g \in E(X_p)$ such that $g_{|_{K}}=k(f)_{|_{K}}$,
$g(x_i)=f(x_i)$ and $g(x) \leq k(f)(x)- \frac{\varepsilon}{2}$ for some $x \in X_0 \setminus K$.\\
Let $L_0$ be a compact set witnessing the fact that $k(f)$ is
$2\varepsilon_0$-saturated, and choose $z_0 \in L_0$; there exists
$f_1 \in E(X_p)$ such that ${f_1}_{|_{L_0}}=k(f)_{|_{L_0}}$,
$f_1(x_i)=f(x_i)$ for $i=1 \ldots n$ and $z_1 \in X_0 \setminus
L_0$ such that
$f_1(z_1) \leq \min \{k(f) (z)+d(z,z_1)\colon z \in L_0\} - \frac{\varepsilon_0}{2} $. \\
Again, let $\varepsilon_1=\inf\{\varepsilon>0 \colon
{f_1}_{|_{X_0}} \mbox { is } \varepsilon-\mbox{saturated } \}$: if
$\varepsilon_1=0$ we are finished, so assume it is not, let $L_1
\supseteq B(z_0,\mbox{diam}(L_0)+d(z_0,z_1))$ be a compact set
witnessing the fact that $f_1$ is $2\varepsilon_1$-saturated and
apply the same process as above to $(f_1,L_1,\varepsilon_1)$.\\
\noindent Then we obtain $z_2 \not \in L_1$ and $f_2 \in E(X_p)$
such that $f_2(x_i)=f(x_i)$ for $i=1 \ldots n$,
${f_2}_{|_{L_1}}={f_1}_{|_{L_1}}$ and $f_2(z_2) \leq \min
\{f_1(z)+d(z,z_2)\colon z \in L_1\} - \frac{\varepsilon_1}{2} $.
\\

\noindent We may iterate this process, thus producing a (finite or
infinite) sequence $(f_m) \in E(X_p)$ who has (among others) the
property that \\
$f_m(x_i)=f(x_i)$ for all $m$ and $i=1 \ldots n$;
the process terminates in finite time only if some
${f_m}_{|_{X_0}}$ is saturated, in which case we have won.\\
So we may focus on the case where the sequence is infinite: then
the construction produces a sequence of $\varepsilon_m$- saturated
Kat\v{e}tov maps $(f_m)$, an increasing and exhaustive sequence of
compact sets $(L_m)$ witnessing that $f_m$ is
$2\varepsilon_m$-saturated, and points $z_m \in L_m \setminus
L_{m-1}$ such that \\
$\ds{f_m(z_m) \leq \min \{ f_{m-1}(z) + d(z,z_m)
\colon z \in L_{m-1} \} - \frac{\varepsilon_{m-1}}{2} }$.\\
If $0$ is a cluster point of $(\varepsilon_n)$, passing to a
subsequence if necessary, we may apply lemma \ref{saturated}(3)
and thus obtain a map $h \in E(X_p)$
such that $h(x_i)=f(x_i)$ for al $i=1 \ldots n$ and $h_{|_{X_0}}$ is saturated. \\

\noindent Therefore, we only need to deal with the case when there
exists $\alpha >0$ such that $\varepsilon_n \geq 2 \alpha$ for all
$n$; we will show by contradiction that this never happens.\\
Since the sequence $(L_m)$ is exhaustive, $(f_n)$ converges
pointwise to some $h \in E(X_p)$ such that
$h(z_m)=f_m(z_m)$ for all $m$.  \\
Up to some extraction, we may assume, since $X$ is compactly
tentacular,
that for all $m$ we have \\
$d(z_0,z_m)+d(z_m,z_{m+1}) \leq d(z_0,z_{m+1})+\frac{\alpha}{2}$.\\
Also we know that $h(z_{m+1})\leq h(z_m)+d(z_m,z_{m+1})-\alpha$.
\\
The two inequalities combined show that $h(z_{m+1})-d(z_{m+1},z_0)
\leq h(z_m)-d(z_m,z_0)-\frac{\alpha}{2}$. \\
This is clearly absurd, since if it were true the sequence
$(h(z_m)-d(z_m,z_0))$ would have to be unbounded, whereas we
have necessarily $h(z_m)-d(z_m,z_0) \geq -h(z_0)$.\\
This is enough to  conclude the proof. $\hfill \lozenge$\\

\noindent {\bf Remark.}  If one applies the construction above to
$X_0=(\N,|.|)$, one obtains a countable set $\{x_n\}_{n \in \N}
\subseteq \U$ such that $d(x_n,x_m)=|n-m|$ for all $n, \ m$ and
$$\forall z \in \U \, \forall \varepsilon >0 \, \exists n,m \in \N \ d(x_n,z)+d(z,x_m) \leq |n-m|+\varepsilon .$$
In particular, $\{x_n\}$ is an isometric copy of $\N$ which is not
contained in any isometric copy of $\R$. \\

\end{subsection}
\end{section}

$ $\\

\noindent  Equipe d'Analyse Fonctionelle, Universit\'e Paris 6 \\
Bo\^ite 186, 4 Place Jussieu, Paris Cedex 05, France.\\
e-mail: melleray@math.jussieu.fr
\end{document}